\documentclass[10pt,a4paper]{article}
\usepackage[latin1]{inputenc}
\usepackage{amsmath}
\usepackage{amsfonts}
\usepackage{amssymb}
\usepackage{graphicx}
\author{Uboho Unyah}
\date{}
\title{An Extension Of Vinogradov's Theorem\footnotetext{This work is dedicated to Aniefiok and Nkoi (Nse), those brothers who got missing from this world in May and July 2007. And from then our feasting days were no longer merry, our holidays were no longer sacred. They laid my foundation before school days, they in their scientific experiments made me see science as an everyday experience. His painting: "INSPIRATION" is still here with us,  inspiring us everyday.}}
\begin{document}
	\maketitle
	\begin{abstract}
		In 1937 Ivan Vinogradov proved the three prime sum version of the Goldbach Conjecture\cite{tA76}, often called the weak form of Goldbach Conjecture. And that it holds for "sufficiently large" odd natural numbers. In this work we use Dirichlet Theorem, Modulo Arithmetics, etc. to extended Vinogradov's Theorem such that every sufficiently large natural number (both even and odd) can be expressed as a sum of three primes. We highlight the configuration of primes for any special case of the three prime sum.   Hence we obtain  Vinogradov's Theorem as a special case of this extended version. We show how Vinogradov's Theorem implies the Goldbach Conjecture; how it (Vinogradov's Theorem) can be derived from it and vice versa. We also obtain the lower bound of the sufficiently largeness. And concluding, we highlight some relationships between the partition function of the Vinogradov's integer $w(v)$ and the Goldbach integer $w(m)$. 
	\end{abstract}
\section{Introduction}
Vinogradov's Theorem as popularly known states that every odd $v\in\mathbb{N}$, $v$ sufficiently large can be expressed as $v=p_1+p_2+p_3$ such that $p_1,p_2,p_3\in\mathbb{P}$. This theorem is  widely known and has been instrumental to numerous mathematical advances\cite{BT04},\cite{mN96},\cite{mN00}. On the major, these important results were obtained using analytical methods. In this work we extend the Vinogradov's Theorem using an elementary method based on the Dirichlet's Theorem and the Modulo Arithmetic. Here our extension of the named theorem follows similarly from the methods we adopted for our prove of Goldbach Conjecture.\cite{uU21c}.\\
\\In this work we obtain an important result that\begin{equation}
p_1+p_2+p_3=2n+b
\end{equation} with $n,b\in\mathbb{N}$, where $b$ satisfies certain prescribed conditions determined by the congruence of $p_3\pmod2$ and $n$ is bounded below based on a condition determined by the set $\mathbb{P}$. And also we obtain a corollary that:
\begin{equation}
m+p_3=2n+b
\end{equation} by making use of our proof of the Goldbach Conjecture itself. So that $m$ in this case is an even integer that can be expressed as a sum of two primes. (The primes here satisfying a number of related conditions; one of such being the Dirichlet's Theorem.)\\
\\We will find Dirichlet's Theorem quite indispensable, since it gives us the basis to prove as in \cite{uU21c} that if $p_i\equiv p_r\pmod2$ then $p_i-p_r=2k$ and $p_i+p_r=2n$ ($n\neq k$). This gives us the corollary that $v=m+p_3$. Using this corollary, we establish the partition function of $v$, $w(v)$ as it relates to (and as a function of) $w(m)$, the partition function of $m$.
\section{Primes And Congruences}
To construct an extension of the weak form of the Goldbach Conjecture (i.e. the three prime sum version of the Goldbach Conjecture) as is defined by Vinogradov's Theorem, we will in this case apply a similar approach as we did in our prove of the strong form of it. Here the concept of primes and some properties of the modulo arithmetic will be quite indispensable. Now let's look at some important definitions and then proceed from there. 
\subsection{Definition}
Let $p\in \mathbb{N}$ such that $p>1$, $p$ is said to be prime if its only  divisors are $1$ and $p$\cite{mN00}. The set of primes is  denoted $\mathbb{P}=\{2,3,5,7,....\}$.
\subsection{Definition}
Let $a,b\in\mathbb{Z}$. The integers $a$ and $b$ are said to be congruent modulo $c$ denoted by $a\equiv b\pmod c$, if $c$ divides $(a-b)$.\\
\\Now let's $P$ be a subset of primes, such that $P=\{p|p<m,p\in\mathbb{P} \}$ for some $m\in\mathbb{N}$, so that $P$ contains all primes less than $m$. Recalling that in $P$ only one element is even, we can choose to reconstruct $P$ so that \begin{equation}
P=P\backslash\{p_k\}\cup\{p_k\}
\end{equation}
where $\{p_k\}$ is the only even prime in this case. Of course we already know that $P\backslash\{p_k\}\cap\{p_k\}=\emptyset$ and as well, we can define $V$ to be \begin{equation}
V=(P\backslash\{p_k\})^2\cup(\{p_k\})^2
\end{equation}
so that \begin{eqnarray}
V&=&(P\backslash\{p_k\})^2\cup(\{p_k\})^2
\\&=&\{(p_i,p_r)|(p_i,p_r)\in(P\backslash\{p_k\})^2\vee(p_i,p_r)\in(\{p_k\})^2 \}
\end{eqnarray}
with $p_k$ the only even prime as before, we can rewrite $V$ as \begin{eqnarray}
V&=&\{(p_i,p_r)|p_i,p_r\in[1]_2\vee p_i,p_r\in[0]_2 \}
\\&=&\{(p_i,p_r)|p_i\equiv p_r\pmod2 \}
\end{eqnarray}
here $[1]_2,[0]_2$ are the residue classes modulo $2$ respectively.\\
\\That is, the prime pair $(p_i,p_r)$ is in $V$ if and only if they are congruent  $\pmod2$. It is important to note that $p_i,p_r\in\{p_k\}$ if and only if $p_i=p_r=p_k$. We find the set $V$ an important asset in constructing our proof of the Goldbach Conjecture. This was treated elaborately in \cite{uU21c}.  At this  point we can take a look at an important theorem that makes an effective use of the concept of  modulo arithmetic; the Dirichlet's Theorem on primes in arithmetic progression.
\subsection{Dirichlet's Theorem}
Let $a,d\in\mathbb{Z}$ with $(a,d)=1$ then there are infinitely many primes satisfying the congruence $p_i\equiv a \pmod d$.\\
Proof:\\
The proof of this theorem can be found in \cite{dE66},\cite{tA76}.$\square$
\subsection{Proposition (Maillet's Conjecture)}\label{P:2.4} 
 There  exist at least a pair of primes $p_i,p_r\in A$ that differ by $h$, where  $h$ is any multiple of $2$. So that  $p_i-p_r=h$; $h=2k$, $\forall k\in\mathbb{Z}$.\cite{uU21c},\cite{jP12}\\
Proof:\\
Let $p_i\equiv a \pmod d$ and $p_r\equiv a \pmod d$ be two primes satisfying Dirichlet's Theorem.\\
 Then \begin{eqnarray}
 p_i\equiv p_j&\equiv& a \pmod d
 \\\Rightarrow p_i&\equiv& p_r \pmod d
 \\\Rightarrow p_i-p_r&\equiv&0  \pmod d
 \end{eqnarray}  
with $d=2$, we have that 
\begin{eqnarray}
p_i-p_r&\equiv& 0 \pmod 2
\\\Rightarrow p_i-p_r&=&2k, \forall k\in\mathbb{Z}
\end{eqnarray}
recall that $(p_i,p_r,d)=1$ so that defining $Z(x,y):=x-y$, we have that there exist at least a prime pair $p_i,p_r\in\mathbb{P}$ such that $Z(p_i,p_r)=2k$ for all $k\in\mathbb{Z}$.$\square$\\
\\$\boldsymbol{Note:}$ It is readily observable that if $Z(x,y)=x-y=2k$, we know that $Z(y,x)=y-x=-(x-y)=-2k$ and that $Z(-x,-y)=Z(x,y)$; so that $Z(-x,-y)=Z(x,y)\Leftrightarrow(-x,-y)\sim_z(x,y)$. Also if $(p_i,p_r),...,(p_s,p_j)$ are ordered pairs with $Z(p_i,p_r)=...=Z(p_s,p_j)$ such that $(p_i,p_r)\neq...\neq(p_s,p_j)$ since each ordered pair is unique. Then each pair generating the same integer is said to be equivalent. As we have already treated concept of prime pair equivalences in \cite{uU21}. Hence we will not dwell on it at this point.\\
\\ Now let's take a look at another proposition that is qualitatively similar into proposition \eqref{P:2.4}.  
\subsection{Proposition}\label{2.5}
 	There  exist at least a pair of primes $p_i,p_j \in A$ that sum to $m$, where $m$ is any multiple of $2$. So that $p_i+p_j=m$; $m=2n$,  $\forall n\in\mathbb{Z}$.\cite{uU21c}\\
Proof:\\
Let $p_i\equiv a \pmod d$ and $p_j\equiv a \pmod d$ be two primes satisfying Dirichlet's Theorem.\\
Then with $d=2$, \begin{eqnarray}
p_i\equiv p_j&\equiv& a \pmod 2
\\p_i&\equiv& p_j\pmod 2
\\p_i+p_j\equiv p_j+p_j\equiv 2p_j &\equiv&0 \pmod 2
\end{eqnarray}
(since $2x\equiv 0\pmod2$ for every $x\in\mathbb{Z}$; $p_j\in\mathbb{P}\subset\mathbb{Z}$)\\
and as before $(p_i,p_j,d)=1$ so that\begin{eqnarray}
p_i+p_j&\equiv& 0\pmod2
\\p_i+p_j&=&2n,\forall n\in\mathbb{Z}
\end{eqnarray}
As required by Goldbach Conjecture, $n\geq2$ since $\forall p\in\mathbb{P}$, $p\geq2$. This forms the central point of our work in \cite{uU21c}. And from here, it is clear that if $n<2$ then at least one prime in the pair will be negative. But with an extension to the set $\mathbb{P}^*$ such that $\mathbb{P}^*=\{p,-p|p\in\mathbb{P}\}$ contains all primes (both positive and negative), then the equation\begin{equation}
p_i+p_j=2n
\end{equation} 
holds for all $n\in\mathbb{Z}$ as required. \\
\\At this point we will adapt these propositions and  employ them in extending the Vinogradov's Theorem; the three prime version of Goldbach Conjecture. 
\section{The Sums Of Three Primes}
In this section we will state and prove the most important result of this work; a theorem on the sum of three primes. This as we had earlier discussed, forms an extension of the Vinogradov's Theorem. We will then take a brief look at some properties of this result and some inherent restrictions as well. 
\subsection{Theorem}\label{3.1}
Every sufficiently large natural number (both and odd) is the sum of three primes.\\
Proof:\\
Let $p_1,p_2$ be primes satisfying Dirichlet's Theorem. So that $p_1\equiv a \pmod d$ and $p_2\equiv a \pmod d$. Let $a,d\in\mathbb{Z}$, $(a,d)=1$ where $d=2$. Let $p_3$ be another prime (not necessarily distinct) such that $p_3\equiv b \pmod d$ where \begin{equation}
b =
\left\{ \begin{array}{ll}
1 & \mbox{if $p_3$ is odd} \\
0 & \mbox{if $p_3$ is even}
\end{array}
\right.
\label{eqtype2}
\end{equation}
Then from the equations,
\begin{eqnarray}
p_1\equiv a\pmod2
\\p_2\equiv a\pmod2
\\p_3\equiv b\pmod2
\end{eqnarray}
Adding up, we obtain \begin{eqnarray}
p_1+p_2+p_3\equiv 2a+b\equiv b \pmod2
\end{eqnarray}
(since $2a\equiv 0\pmod2$ for all $a\in\mathbb{Z}$)\\
And \begin{eqnarray}
p_1+p_2+p_3&\equiv& b\pmod2
\\p_1+p_2+p_3-b&\equiv& 0 \pmod2
\\p_1+p_2+p_3-b&=&2n, \forall n\in\mathbb{Z}
\\p_1+p_2+p_3&=&2n+b, \forall n\in\mathbb{Z}
\end{eqnarray}
Now with $p_1,p_2,p_3\in\mathbb{P}$ we of course know that $\forall p\in\mathbb{P}$, $p\geq2$. So that \begin{eqnarray}
p_1+p_2+p_3&=&2n+b, \forall n\in\mathbb{Z}
\\\Rightarrow p_1+p_2+p_3&=&2n;\, (b=0 \, \text{since }\,  p_3\,  \text{is even})
\\\Rightarrow p_1+p_2+p_3&=&2n\geq 2+2+2=6
\\\Rightarrow n\geq3, n\in\mathbb{Z}
\end{eqnarray}
Hence we derive that \begin{equation}
p_1+p_2+p_3=2n+b, b =
\left\{ \begin{array}{ll}
1 & \mbox{if $p_3$ is odd} \\
0 & \mbox{if $p_3$ is even}
\end{array}
\right.\forall n\in\mathbb{Z},n\geq3 
\label{eqtype1}
\end{equation}
That is, every sufficiently large natural number (odd or even) can be expressed as a sum of three primes.$\square$\\
\\We can verify the proof as we have just seen. Using the insights we have gained from a proof of Goldbach Conjecture which we gave recently in \cite{uU21c}. We can recall from proposition \eqref{2.5} that there exist at least a pair of primes $p_i,p_j$ such that $p_i+p_j=2n$, $\forall n\in\mathbb{Z}$; with the condition that $p_i\equiv p_j\pmod2$, so that the Dirichlet's Theorem's sufficient condition is satisfied. That is, $p_i\equiv a\pmod2$, $p_j\equiv a\pmod2$ with $(a,2)=1$.\\
\\For the Goldbach Conjecture, the restriction that $n\geq2$ is that $p_i,p_j\in\mathbb{P}$ (where $\forall p\in\mathbb{P}$, $p\geq2$) so that $p_i+p_j\geq4$. Now using Goldbach Conjecture as a guide, we observe that \begin{equation}
p_1+p_2+p_3=2n+b, b =
\left\{ \begin{array}{ll}
1 & \mbox{if $p_3$ is odd} \\
0 & \mbox{if $p_3$ is even}
\end{array}
\right.\forall n\in\mathbb{Z},n\geq3 
\label{eqtype5}
\end{equation} can be written as \begin{equation}
m+p_3=2n+b, b =
\left\{ \begin{array}{ll}
1 & \mbox{if $p_3$ is odd} \\
0 & \mbox{if $p_3$ is even}
\end{array}
\right.\forall n\in\mathbb{Z},n\geq3 
\label{eqtype6}
\end{equation} 
Here $m=p_1+p_2$ as in Goldbach Conjecture. We can also recall that $p_1\equiv a\pmod2$ and $p_2\equiv a\pmod2$, so that
\begin{equation}
p_1\equiv a\equiv p_2\pmod2\Rightarrow p_1\equiv p_2\pmod2
\end{equation} 
with $(a,2)=1$ etc. as shown before. This equation \eqref{eqtype6} quite obviously forms a corollary for every sufficiently large integer. For an example, every sufficiently large even integer can be expressed as $v=m+2$ and every sufficiently large odd integer can be expressed as $v=m+3$ etc.\\
\\$\boldsymbol{Remark:}$ We can easily derive Vinogradov's Theorem from here, by restricting $p_3$ to only odd primes. So that   \begin{equation}
p_1+p_2+p_3=2n+1, \forall n\in\mathbb{Z}, n\geq3,p_3\geq3
\end{equation}
where the two primes $p_1,p_2$ need only to be congruent $\pmod2$, that is $p_1\equiv a\equiv p_2\pmod2$ with $(a,2)=1$ as required by Dirichlet's Theorem.\\
The three prime sum for even integers follows as before except that in this case, $p_3$ is strictly even so that $b=0$ (since $p_3\equiv b\pmod2$), giving\begin{equation}
p_1+p_2+p_3=2n, \forall n\in\mathbb{Z}, n\geq3,p_3=2
\end{equation} 

\subsection{Goldbach-Vinogradov Relationships}
Vinogradov's Theorem has always been known to imply the Goldbch Conjecture. Given an integer large enough,  
\begin{equation}
p_1+p_2+p_3=2n+b, b =
\left\{ \begin{array}{ll}
1 & \mbox{if $p_3$ is odd} \\
0 & \mbox{if $p_3$ is even}
\end{array}
\right.\forall n\in\mathbb{Z},n\geq3 
\label{eqtype4}
\end{equation}
For $p_3$ even, $p_3=2k=2$ (since only $p=2$ is even prime); so that \begin{equation}
p_1+p_2+2=2n\Rightarrow p_1+p_2=2n-2=2(n-1)
\end{equation}
with $2(n-1)=m$ we have $p_1+p_2=m$ as usual for Goldbach Conjecture.\\ Also for $p_3$ odd, $p_3=2k+1$ for some $k$; so that $b=1$ then \begin{eqnarray}
p_1+p_2+p_3&=&2n+1
\\\Rightarrow p_1+p_2&=&2n+1-p_3
\\&=&2n+1-(2k+1)
\\\Rightarrow p_1+p_2&=&2n-2k
\\&=&2(n-k)
\end{eqnarray}
and with $2(n-k)=m$ we obtain 
\begin{equation}
p_1+p_2=2(n-k)=m
\end{equation}
as in Goldbach Conjecture.\\
\\

We can infer that from \begin{equation}
p_1+p_2+p_3=2n+b,b =
\left\{ \begin{array}{ll}
1 & \mbox{if $p_3$ is odd} \\
0 & \mbox{if $p_3$ is even}
\end{array}
\right.\forall n\in\mathbb{Z},n\geq3 
\label{eqtype3}
\end{equation}
we have \begin{equation}
m+p_3=2n+b
\end{equation}
for all sufficiently large $n$. That is, we can derive as a corollary that  every sufficiently large natural number (both even and odd), is the sum of an even integer and a prime. This is true even if we fix $p_3$ to be some specific value  of prime, we can still draw the same conclusions. We can take a look at some specific examples:\begin{itemize}
	\item for $p_3=2$: $m+2=2n$ for all sufficiently large $n$.\\
	\subsubsection{Example:} \subitem $4+2=6$
	\subitem $6+2=8$
	\subitem $8+2=10$
	\subitem $10+2=12$
	\subitem $....$
	\subitem $....$
	\item for $p_3=3$: $m+3=2n+1$ for all sufficiently large $n$.\\
	\subsubsection{Example:} \subitem $4+3=7$
	\subitem $6+3=9$
	\subitem $8+3=11$
	\subitem $10+3=13$
	\subitem $12+3=15$
	\subitem $....$
	\subitem $....$
	\item for $p_3=5$: $m+5=2n+1$ for all sufficiently large $n$.\\
	\subsubsection{Example:} \subitem $4+5=9$
	\subitem $6+5=11$
	\subitem $8+5=13$
	\subitem $10+5=15$
	\subitem $....$
	\subitem $....$
	\item for $p_3=7$: $m+7=2n+1$ for all sufficiently large $n$.\\
	\subsubsection{Example:} \subitem $4+7=11$
	\subitem $6+7=13$
	\subitem $8+7=15$
	\subitem $10+7=17$
	\subitem $....$
	\subitem $....$
\end{itemize}
 This list goes on endlessly. And as well, it is interesting to observe that a sufficiently large $v$ can in itself be expressed as a sum $v=m+p$ in many different ways. 
 \subsubsection{Example:}
 Let $v=21$, we can write this as a sum of an even integer and a prime in many different ways. So that   \begin{eqnarray}
 	21&=&18+3
 	\\&=&16+5
 	\\&=&14+7
 	\\&=&10+11
 	\\&=&8+13
 	\\&=&4+17
 \end{eqnarray}
 This sort of multiplicity holds only for odd $v$, since $p_3$ can take several distinct values which are all odd as required.   
 \section{Partition Function $w(v)$ on $v$}
 Using the insights we had obtained from the partition function $w(m)$ of the Goldbach even integer $m$, we can in this case make some interesting observations about the prime partitions of the Vinogradov integer $v$. We will follow this in two cases: when $v$ is even and when $v$ is odd. But how do we differentiate these cases? First we can recall that as a corollary, every integer that is a sum of three primes can be written as a sum of an even integer and a prime. So that if $v=p_1+p_2+p_3$ then $v=m+p_3$; where $p_3$ is even or odd accordingly.   
 \subsection{When $v$ Is Even}
 Let $v$ be an even integer that can be expressed as a sum of three primes, then from the condition on $b$ as in theorem \eqref{3.1} we can express $v=m+2$, where $p_3=2$ and $p_3=2\equiv b\pmod2$ with $b=0$. \\
 Now with $w(v)$  as the number of expressions of $v$ as a sum of three primes, we have that 
 \begin{equation}
 w(v)=w(m+2)=w(m)
 \end{equation}
 That is, the number of expressions for $v=p_1+p_2+p_3$ is equal to the number of expressions for $m=p_1+p_2$ if $v$ is even. This holds since $p_3=2$ is the only way of expressing $p_3$, given that $p_1,p_2$ by definition must be congruent.
 \subsubsection{Example:}
 Let $v=20$, then we can write it as \begin{equation}
 v=20=18+2
 \end{equation}
 where $m=18$. Now the different expressions for $v=20$ as a sum of three primes can be given as:\begin{eqnarray}
 v=20&=&5+13+2
 \\&=&7+11+2
 \end{eqnarray} 
 Observably we can express $m=18$ as a sum of two primes 
 \begin{eqnarray}
 m=18&=&5+13
 \\&=&7+11
 \end{eqnarray}
 So that if $v=p_1+p_2+p_3$ and $m=p_1+p_2$ with $p_3=2$, then $w(v)=w(m)$ as required. 
 \subsection{When $v$ Is Odd}
 As before, let $v=p_1+p_2+p_3$ be odd, so that $v=m+p_3$ with $p_3$ odd in this case and $b=1$. We know that $p_3\geq3$ so that for sufficiently large $v$, $p_3$ can take different possible values as shown in the examples in the previous section. That is here $v$ can be expressed as:\begin{eqnarray}
 v=m_1+p_1
 \\=m_2+p_2
 \\=m_3+p_3
 \\...........
 \\...........
 \\=m_k+p_k
 \end{eqnarray}
 with $m_i\neq m_j$ and $p_i\neq p_j$.\\
 \\Hence it is obvious that $w(v)\geq w(m)$ for any $m$ such that $v=m+p$. In essence, $w(v)$ in the odd cases is simply the sum of the different $w(m_i)$, $i=1$ to $k$. So that \begin{eqnarray}
 w(v)&=&w(m_1)+w(m_2)+...+w(m_k)
 \\&=&\sum_{i=1}^{k}w(m_i)
 \end{eqnarray} 
 and $w(v)\geq w(m_i)$ for any $i$ and since it can take many different odd values as required. 
 \section{Conclusion}
 The ability to write even integers as a sum of three primes is quite obvious since every even integer can be obtained from the previous by a successive addition of $2$. Yet its quite interesting to realize that one may be able to write an even integer both as a sum of two primes and a sum of three integers as well. An application of this property is found in a more generalized form where an integer sufficiently large can be written as the sum of as many primes as possible. We will look at this generalized form in \cite{uU21d}. But for the present work as we have succeeded to show, writing every sufficiently larger natural number (both even and odd) as a sum of three primes pose no contradiction.  
 
 Department Of Mathematics, University Of Uyo, Nigeria\\
 \text{Email:unyahinmaths2020@gmail.com}


\begin{thebibliography}{99}
 	\bibitem{tA76}
 	Apostol, Tom: {Introduction To Analytic Number Theory}, Springer-Verlag New York, 1976.
 	\bibitem{BT04}
 	Green, Ben. and Tao, Terrence: The Primes Contain Arbitrarily Long Arithmetic Progressions, (Preprint). Available at http://www.arXiv.org
 \bibitem{mN00}
 Nathanson, M. B.: Elementary Methods In Number Theory, Springer-Verlag, New York Inc., New York, 2000.
 \bibitem{mN96}
 Nathanson, M. B.: Additive Number Theory: The Classical Bases, Springer-Verlag, New York Inc., New York, 1996.
 \bibitem{jP12}
 \text{J\'{a}nos} Pintz: On Differences Of Primes,
 http://arxiv.org/abs/1206.0149v1
\bibitem{uU21}
Unyah, Uboho: A Study of Some Equivalence Properties Of Primes (In Their Pairs), arXiv:2110.07334v1[math.GM], 12 Oct 2021.
\bibitem{uU21c}
Unyah, Uboho: A Look At Goldbach Conjecture, arXiv:2110.07334v4[math.GM]. 
\bibitem{uU21d}
Unyah, Uboho:A Generalized Partition Of Integers In Primes. (Preprint)
\bibitem{dE66}
Wang, Zijiang: {Elementary Proof of Dirichlet theorem}. Available at http://math.uchicago.edu/may/REU2017/REUPapers/WangZijian.pdf
 \end{thebibliography}
\end{document}